\documentclass[11pt,twoside]{article}

\usepackage{a4wide}
\usepackage{amsfonts}
\usepackage{amssymb}
\usepackage{amsmath}
\usepackage{graphicx}
\usepackage{nopageno}
\usepackage{xcolor}

\newcommand{\ignore}[1]{}

\def\@begintheorem#1#2{\par\bgroup{\sc #1\ #2. }\it\ignorespaces}
\def\@opargbegintheorem#1#2#3{\par\bgroup{\sc #1\ #2\ (#3). } \it\ignorespaces}
\def\@endtheorem{\egroup}
\newtheorem{theorem}{Theorem}[section]
\newtheorem{corollary}[theorem]{Corollary}
\newtheorem{conjecture}[theorem]{Conjecture}
\newtheorem{question}[theorem]{Question}
\newtheorem{lemma}[theorem]{Lemma}

\newtheorem{example}[theorem]{Example}
\newtheorem{proposition}[theorem]{Proposition}
\newtheorem{definition}[theorem]{Definition}
\newcommand{\bt}[1]{\begin{theorem}\label{#1}}
\newcommand{\bc}[1]{\begin{corollary}\label{#1}}
\newcommand{\bconj}[1]{\begin{conjecture}\label{#1}}
\newcommand{\bq}[1]{\begin{question}\label{#1}}
\newcommand{\bl}[1]{\begin{lemma}\label{#1}}
\newcommand{\be}[1]{\begin{example}\label{#1}}
\newcommand{\bp}[1]{\begin{proposition}\label{#1}}
\newcommand{\ba}[1]{\begin{algorithm}\rm\label{#1}}
\newcommand{\bd}[1]{\begin{definition}\rm\label{#1}}{\normalsize }
\newcommand{\bpr}{\noindent {\em Proof. }}
\newcommand{\et}{\end{theorem}}
\newcommand{\ec}{\end{corollary}}
\newcommand{\econj}{\end{conjecture}}
\newcommand{\eq}{\end{question}}
\newcommand{\el}{\end{lemma}}
\newcommand{\ee}{\end{example}}
\newcommand{\ep}{\end{proposition}}
\newcommand{\ed}{\end{definition}}
\newcommand{\epr}{{\ \vbox{\hrule\hbox{%
\vrule height1.3ex\hskip0.8ex\vrule}\hrule}}\\\par}
\newcommand{\mepr}{{\ \ \ \vbox{\hrule\hbox{%
\vrule height1.3ex\hskip0.8ex\vrule}\hrule}}}

\def\R{\mathbb{R}}
\def\H{{\cal H}}
\def\G{{\cal G}}
\def\HC{{\overline\H}}

\begin{document}

\title{\bf Separable and Equatable Hypergraphs}

\author{
Daniel Deza
\thanks{\small University of Toronto. Email: daniel.deza@mail.utoronto.ca}
\and
Shmuel Onn
\thanks{\small Technion - Israel Institute of Technology.
Email: onn@ie.technion.ac.il}
}
\date{}

\maketitle

\begin{abstract}

We consider the class of {\em separable} $k$-hypergraphs, which can be viewed as uniform
analogs of threshold Boolean functions, and the class of {\em equatable} $k$-hypergraphs.
We show that every $k$-hypergraph is either separable or equatable but not both.
We raise several questions asking which classes of equatable (and separable) hypergraphs
enjoy certain appealing characterizing properties, which can be viewed as uniform
analogs of the $2$-summable and $2$-monotone Boolean function properties. In particular,
we introduce the property of {\em exchangeability}, and show that all these questioned
characterizations hold for graphs, multipartite $k$-hypergraphs for all $k$, paving $k$-matroids
and binary $k$-matroids for all $k$, and $3$-matroids, which are all equatable if and only if they
are exchangeable. We also discuss the complexity of deciding if a hypergraph is separable, and in
particular, show that it requires exponential time for paving matroids presented by independence
oracles, and can be done in polynomial time for binary matroids presented by such oracles.

\end{abstract}

\section{Introduction}

A {\em $k$-hypergraph} is a pair $(V,\H)$ where $V$ is a finite set and $\H$ is a set of
$k$-subsets of $V$. The elements $v\in V$ are {\em vertices} and the sets $E\in\H$
are {\em edges}. We avoid unnecessary trivialities and in all our statements,
questions, lemmas, and theorems, assume $1\leq k<n:=|V|$.

Here we study the following two classes of hypergraphs. First, a $k$-hypergraph is
{\em separable} if there is a labeling $x:V\rightarrow\R$ of vertices by real numbers such that,
with $x(E):=\sum_{v\in E}x(v)$,
\begin{equation}\label{separable_definition}
\H\ =\ \{E\subseteq V\ :\ |E|=k,\ \ x(E)\geq 0\}\ ,
\end{equation}
i.e., there is a vertex labeling $x$ such that the edges of
$\H$ are the $k$-sets $E$ with nonnegative sum $x(E)$.
For instance, the following $3$-hypergraph $(V,\H)$ is separable with a suitable labeling $x$,
$$V=[6],\quad\H=\{124,134,145,146,234,245,246,345,346\},\quad
x(v)=
\left\{
  \begin{array}{rl}
    -1, & v=1,2,3; \\
    3, & v=4; \\
    -2, & v=5,6.
  \end{array}
\right.
$$
Second, a $k$-hypergraph is {\em equatable} if there is a nonzero labeling
$y:{V\choose k}\rightarrow\R_+$ of $k$-subsets of $V$ by nonnegative real numbers such that,
with $\HC:={V\choose k}\setminus\H$ the {\em complement} of $\H$,
\begin{equation}\label{equatable_definition}
\mbox{for every $v\in V$ there holds}
\ \ \sum\{y(E)\ :\ v\in E\in \H\}\ =\ \sum\{y(F)\ :\ v\in F\in \HC\}\ ,
\end{equation}
i.e., there is a labeling of all $k$-subsets such that, for every vertex $v$,
the sum of labels of the edges $E\in\H$ containing $v$ is equal to the sum of labels
of the non edges $E\in\HC$ containing $v$. For instance, the following $3$-hypergraph
$(V,\H)$ is equatable with a suitable labeling $y$,
$$V=[6],\quad\H=\{135,136,145,146,235,236,245,246\},\quad
y(G)=
\left\{
  \begin{array}{rl}
    1, & G=134,135,246,256; \\
    0, & \hbox{otherwise}.
  \end{array}
\right.
$$
As we show in Lemma \ref{separable_equatable}, the Farkas' lemma in linear programming implies that
every $k$-hypergraph is either seprabale or equatable but not both. So equatable $k$-hypergraphs
are exactly the non separable ones, but it is convenient to have both terms at hand.
As we show in Lemma \ref{sufficient}, if a $k$-hypergraph $(V,\H)$ has edges $E_1,E_2\in\H$ and
non edges $F_1,F_2\in\HC$ such that $E_1\cap E_2=F_1\cap F_2$ and $E_1\cup E_2=F_1\cup F_2$,
then the hypergraph is equatable. Thus, it is natural to ask in which classes, equatable hypergraphs
are characterized by this property.

\bq{question1}
In which classes, a $k$-hypergraph $(V,\H)$ is equatable if and only if there
are edges $E_1,E_2\in\H$ and non edges $F_1,F_2\in\HC$ such that
$E_1\cap E_2=F_1\cap F_2$ and $E_1\cup E_2=F_1\cup F_2$?
\eq

Separable hypergraphs are related to {\em threshold Boolean functions} which had been extensively
studied in the literature and have numerous applications, hence the interest in their study,
see \cite{CH2} and the references therein. Boolean functions $f:2^V\rightarrow\{0,1\}$ are in
correspondence with hypergraphs $(V,\H)$ via $\H:=\ker(f)=\{E\subseteq V:f(E)=0\}$.
The function $f$ is {\em threshold} if there is a labeling $x:V\rightarrow\R$
and a real number $t$ with $\ker(f)=\{E\subseteq V:x(E)\geq t\}$.
However, the hypergraphs corresponding to threshold Boolean functions typically contain
edges of different cardinalities, so are {\em non uniform} and not $k$-hypergraphs.
Nonetheless, the non uniform analog of Question \ref{question1} was raised many times in the
Boolean function literature, and asks in which classes, a Boolean function is not threshold if and only
if it is so-called {\em $2$-summable}. This is not true in general, and over the years
many authors had independently discovered examples of functions which are neither threshold nor
$2$-summable. See \cite[Section 9.3]{CH2} for a thorough discussion of these notions and examples.

\vskip.2cm
Here we study the following question with another natural condition which is stronger than the
condition of Question \ref{question1}. Call a $k$-hypergraph $(V,\H)$ {\em exchangeable} if there
are $E_1,E_2\in\H$ and $v_1\in E_1\setminus E_2$, $v_2\in E_2\setminus E_1$ such that
$E_1\setminus\{v_1\}\cup\{v_2\}$, $E_2\setminus\{v_2\}\cup\{v_1\}\in\HC$.

\bq{question2}
In which classes, a $k$-hypergraph is equatable if and only if it is exchangeable?
That is, in which classes, a $k$-hypergraph $(V,\H)$ is equatable if and only if there
are $E_1,E_2\in\H$ and $v_1\in E_1\setminus E_2$, $v_2\in E_2\setminus E_1$ such that
$E_1\setminus\{v_1\}\cup\{v_2\}$, $E_2\setminus\{v_2\}\cup\{v_1\}\in\HC$? \break
(Equivalently, in which classes, a $k$-hypergraph $(V,\H)$ is separable if and only if
for every $E_1,E_2\in\H$ and $v_1\in E_1\setminus E_2$, $v_2\in E_2\setminus E_1$, either
$E_1\setminus\{v_1\}\cup\{v_2\}\in\H$ or $E_2\setminus\{v_2\}\cup\{v_1\}\in\H$?)
\eq

As mentioned, the condition in Question \ref{question2} implies that of Question \ref{question1},
by simply taking $F_1:=E_1\setminus\{v_1\}\cup\{v_2\}$ and $F_2:=E_2\setminus\{v_2\}\cup\{v_1\}$,
but they coincide for $k\leq 3$. It turns out that exchangeability contrasts with a uniform analog
of the so-called {\em monotone} property of Boolean functions. Indeed, in Section 7 we show
that a $k$-hypergraph is {\em not} exchangeable if and only if it is $2$-monotone.
Again, it had been asked over the years, see \cite[Section 8.8]{CH2} and \cite{Gol}, in which
classes is it true that being threshold is equivalent to being $2$-monotone, and several independent
non threshold $2$-monotone examples were given. In Section 6 we give a simple interpretation,
using equatability, of an example from \cite{RRST}, showing an equatable but not exchangeable
$3$-hypergraph, so one violating the conditions of Questions \ref{question1} and \ref{question2}.

\vskip.2cm
Our main result here is the following theorem (see Sections 3--5 for the precise statements
and for the definitions of multipartite hypergraphs, paving matroids, and binary matroids).
\bt{main}
In each of the following classes, equatability is equivalent to exchangeability:
\begin{itemize}
\item
All $k$-hypergraphs for $k=1,2$, and in particular all graphs;
\item
All multipartite $k$-hypergraphs for all $k$;
\item
All paving $k$-matroids and all binary $k$-matroids for all $k$, and all $3$-matroids.
\end{itemize}
\et

\vskip.2cm
The {\em dual} of the $k$-hypergraph $(V,\H)$ is the $(n-k)$-hypergraph $(V,\H^*)$
with set of edges $\H^*:=\{V\setminus E:E\in\H\}$.
As follows from Lemma \ref{complement_dual}, the theorem implies also the following.

\bc{corollary}
Equatability is equivalent to exchangeability also for all complements and all duals of
hypergraphs appearing in Theorem \ref{main}. In particular, this holds for all hypergraphs
with $n=k+1,k+2$ vertices, and all duals of paving $k$-matroids and all duals of $3$-matroids.
\ec

The condition in Question \ref{question1} asserts (see proof of Lemma \ref{sufficient}) that, if a
hypergraph is equatable with a labeling $y$, then it has in fact a labeling with only $0,1$ values
and only $4$ non zeros. While this is not true for all hypergraphs as demonstrated by the example in
Section 6, this example does not exclude a possible positive answer to the following further question.
\bq{question3}
Is it true that $k$-hypergraph $(V,\H)$ is equatable if and only if there
is a labeling $y:{V\choose k}\rightarrow\{0,1\}$ such that \eqref{equatable_definition} holds?
Is there always such $y$ with at most $2k$ nonzero values?
\eq

The property of admitting such a restricted, $0,1$ valued labeling $y$ which satisfies
\eqref{equatable_definition}, is superficially reminiscent of the recent notion of a
{\em null $k$-hypergraph} introduced in \cite{FKPT}, which is a hypergraph admitting a
labeling $y$ satisfying \eqref{equatable_definition} and restricted to the values $y(E)=\pm1$ for
$E\in\H$ and $y(F)=0$ for $F\in\HC$. But such a labeling $y$ has negative values, and indeed,
already for $k=3$, deciding if a $3$-hypergraph is null is NP-complete \cite{FKPT}, while deciding
if a $k$-hypergraph is equatable can be done in polynomial time for any fixed $k$, see Section 8.

The class of separable hypergraphs is geometrically natural and extends the well studied class
of threshold graphs \cite{CH1,MP} from $k=2$ to any positive $k$. It is a rich class,
as already for $k=3$, deciding if a separable $3$-hypergraph has a perfect matching
(a set of edges whose disjoint union is $V$, see \cite{Hax}), is NP-complete \cite{Onn}.
This is in contrast with the well known polynomial time solution of the perfect
matching problem for all graphs. In fact, deciding a perfect matching in a separable
$3$-hypergraph is closely related to the classical NP-complete {\em $3$-partition} problem,
which is precisely to decide if there is a perfect matching in a $3$-hypergraph given by
$\H:=\{E\subseteq V:|E|=3,\ x(E)=0\}$ for some given vertex labeling $x$, see \cite{GJ}.

The class of multipartite hypergraphs appearing in Theorem \ref{main}, extending bipartite graphs,
is also rich, as deciding the existence of a perfect matching in a $3$-partite hypergraph is precisely
the classical NP-complete {\em $3$-dimensional matching} problem, see \cite{GJ} again.

Matroids form a class of hypergraphs which is central to combinatorial optimization \cite{Sch}.
The subclass of {\em paving matroids} is broad and in fact conjectured to contain almost all
matroids \cite{MNWW}, and the subclass of {\em binary matroids} is rich and well studied and
includes in particular all {\em graphic matroids}, see Section 5 and \cite{Oxl}.
The class of $3$-matroids, in its oriented refinement, is already universal
for semi-algebraic varieties in a well defined sense, see \cite{BLSWZ}.

Finally, we note that $k$-hypergraphs are very complicated objects already for $k=3$, and already
quite simple problems over them are hard. For instance, deciding degree sequences is NP-complete
for $3$-hypergraphs \cite{DLMO} while polynomial time doable for graphs. Therefore it would be
particularly interesting to resolve Question \ref{question3} at least for $3$-hypergraphs.

\vskip.2cm
The rest of the article is organized as follows. In Section 2 we give some preparatory lemmas used
throughout. In Sections 3--5 we prove that the condition of Question \ref{question2}, and hence also
of Question \ref{question1}, hold for graphs, multipartite hypergraphs, and paving, binary, and rank
$3$ matroids, respectively. In Section 6 we give a simple presentation, using the notion of equatability
and the assertion of Lemma \ref{separable_equatable} that a hypergraph is equatable if and only
if it is non separable, of an example from \cite{RRST}. In Section 7 we discuss uniform monotonicity
of hypergraphs, show that a hypergraph is not exchangeable if and only if it is $2$-monotone, and raise
the question of whether Theorem \ref{main} could be extended to all matroids, that is, whether any
matroid is equatable if and only if it is exchangeable. Finally, in Section 8, we provide some remarks
on the complexity of deciding if a $k$-hypergraph is separable. We show that it can be done in
polynomial time for each fixed $k$, in polynomial time even for variable $k$ for classes satisfying
the condition of Question \ref{question2} by checking exchangeability and $2$-monotonicity,
requires exponential time for paving matroids presented by independence oracles,
and can be done in polynomial time for binary matroids presented by such oracles.

\section{Preparation}

First we prove that the equatable hypergraphs are precisely the non separable hypergraphs.

\bl{separable_equatable}
Every $k$-hypergraph $(V,\H)$ is either separable or equatable but not both.
\el
\bpr
Define a ${V\choose k}\times V$ matrix $A$, and a vector $b\in\R^{V\choose k}$, by
$$A_{G,v}\ :=\ \left\{
    \begin{array}{rl}
      -1, & v\in G\in\H; \\
      1, & v\in G\in\HC; \\
      0, & v\notin G\in{V\choose k},
    \end{array}
  \right.
\quad\quad
b_G\ :=\ \left\{
    \begin{array}{rl}
      0, & G\in\H; \\
      -1, & G\in\HC. \\
    \end{array}
  \right.
$$
Now, one version of the Farkas' lemma in linear programming (cf. \cite[Corollary 7.1e]{Sch}) asserts
that one and only one of the following two systems of linear inequalities has a solution,
\begin{equation}\label{primal}
Ax\ \leq\ b,\ \ \ x\in\R^V\ ,
\end{equation}
\begin{equation}\label{dual}
y^TA\ =\ 0,\ \ \ y\in\R_+^{V\choose k},\ \ \ y^Tb\ <\ 0\ .
\end{equation}
It is easy to see that $x$ is a solution of \eqref{primal} if and only if it is a
labeling showing $(V,\H)$ is separable, and $y$ is a solution of \eqref{dual}
if and only if it is a labeling showing $(V,\H)$ is equatable.
\epr
The proof of this lemma shows that separability and equatability are essentially dual to
each other in the linear programming sense, hence the interest in their simultaneous study.

Next we show that the conditions of Questions \ref{question1} and \ref{question2}
sufficient for a hypergraph to be equatable. In the non uniform analog notions
in the context of Boolean functions, this corresponds to the fact that $2$-summable
Boolean functions are not threshold \cite[Section 9.3]{CH2}.

\bl{sufficient}
If in a $k$-hypergraph $(V,\H)$ there are $E_1,E_2\in\H$, $F_1,F_2\in\HC$ with
$E_1\cap E_2=F_1\cap F_2$ and $E_1\cup E_2=F_1\cup F_2$, then $(V,\H)$ is equatable,
and moreover, has a $0,1$ labeling with at most $4$ non zeros satisfying
\eqref{equatable_definition}. In particular, if it is exchangeable than it is equatable.
\el
\bpr
Define a labeling $y:{V \choose k}\rightarrow\R$ by $y(E_1)=y(E_2)=y(F_1)=y(F_2)=1$
and $y(E)=0$ for all other $E$. Then this labeling shows that the hypergraph is equatable,
since for any $v\in V$:
\begin{itemize}
\item if $v\notin E_1\cup E_2$ then $\sum\{y(E):v\in E\in\H\}=0=\sum\{y(F):v\in F\in\HC\}$;
\item if $v\in (E_1\cup E_2)\setminus(E_1\cap E_2)$ then
$\sum\{y(E):v\in E\in\H\}=1=\sum\{y(F):v\in F\in\HC\}$;
\item if $v\in E_1\cap E_2$ then $\sum\{y(E):v\in E\in\H\}=2=\sum\{y(F):v\in F\in\HC\}$.
\end{itemize}
Clearly, this labeling is $0,1$ with at most $4$ nonzero values.
And, if the hypergraph is exchangeable, with suitable $E_1,E_2,v_1,v_2$,
take $F_1:=E_1\setminus\{v_1\}\cup\{v_2\}$, $F_2:=E_2\setminus\{v_2\}\cup\{v_1\}$.
\epr

The next lemma together with Theorem \ref{main} implies Corollary \ref{corollary}.
\bl{complement_dual}
The following are equivalent: a $k$-hypergraph $(V,\H)$ is exchangeable; its
complement $(V,\HC)$ is exchangeable; its dual $(V,\H^*)$ is exchangeable.
And the following are equivalent: $(V,\H)$ is equatable;
its complement $(V,\HC)$ is equatable; its dual $(V,\H^*)$ is equatable.
\el
\bpr
For the first statement, if $E_1,E_2\in\H$,
$F_1:=E_1\setminus\{v_1\}\cup\{v_2\},F_2:=E_2\setminus\{v_2\}\cup\{v_1\}\in\HC$,
so that $(V,\H)$ has the desired property, then $F_1,F_2\in\HC$ and
$E_1=F_1\setminus\{v_2\}\cup\{v_1\}$,\break $F_2=E_2\setminus\{v_1\}\cup\{v_2\}\in\H$ so that $(V,\HC)$
also has the desired property, and $E^*_1:=V\setminus E_1$,\break $E^*_2:=V\setminus E_2\in\H^*$ and
$F^*_1:=E^*_1\setminus\{v_2\}\cup\{v_1\},F^*_2:=E^*_2\setminus\{v_1\}\cup\{v_2\}\in\overline{\H^*}$
so that $(V,\H^*)$ also has the desired property. The converses are again proved in the same way.

For the second statement, let $y$ be a labeling of $k$-sets for which $(V,\H)$ satisfies
\eqref{equatable_definition}. Then clearly $(V,\HC)$ also satisfies \eqref{equatable_definition}
with $y$. Now define a labeling $y^*$ of $(n-k)$-sets by $y^*(E):=y(V\setminus E)$.
Note that summing \eqref{equatable_definition} over all $v\in V$ and dividing by $k$ we obtain
$\sum\{y(G):G\in\H\}=\sum\{y(G):G\in\HC\}$. Note also that $\overline{\H^*}=\HC^*$.
Then for each $v\in V$,
\begin{eqnarray*}
  \sum\{y^*(E):v\in E\in\H^*\} &=& \sum\{y(V\setminus E):v\in E\in\H^*\}
  \ =\ \sum\{y(G):v\notin G\in\H\} \\
   &=& \sum\{y(G):G\in\H\}-\sum\{y(G):v\in G\in\H\}\ ,
\end{eqnarray*}
\begin{eqnarray*}
  \sum\{y^*(F):v\in F\in\overline{\H^*}\} &=& \sum\{y(V\setminus F):v\in F\in\HC^*\}
  \ =\ \sum\{y(G):v\notin G\in\HC\} \\
   &=& \sum\{y(G):G\in\HC\}-\sum\{y(G):v\in G\in\HC\}\ ,
\end{eqnarray*}
from which we conclude $\sum\{y^*(E):v\in E\in\H^*\}=\sum\{y^*(F):v\in F\in\overline{\H^*}\}$.
Therefore we find that $(V,\H^*)$ satisfies \eqref{equatable_definition} with the labeling $y^*$.
The converses are proved in the same way.
\epr

In what follows, for a $k$-subset $\{v_1,v_2,\dots,v_k\}$,
we also use the abridged form $v_1v_2\cdots v_k$.

\section{Graphs}

The resolution of Question \ref{question2} for all $k$-hypergraphs with $k=1,2$ is quite simple
and for $k=2$, that is, {\em graphs}, equivalent to various characterizations of {\em threshold graphs}
available in literature. These graphs were introduced in \cite{CH1} and further studied by many authors
and in the books \cite{Gol,MP}. A $2$-hypergraph $(V,\H)$, that is, a graph, is {\em threshold} if
there is a labeling $x:V\rightarrow\R$ and a real number $t$ such that $\H=\{E\subseteq V:x(E)\geq t\}$.
It is easy to see that a graph is threshold if and only if it is separable, that is, has such a
labeling with $t=0$. Before stating the theorem on graphs, we record another well known characterization
of threshold and separable graphs, which will be used later on. A graph $(V,\H)$ is {\em orderable}
if there is an ordering $v_1,v_2,\dots$ of $V$ with each $v_j$ either {\em isolated}, meaning that
$\{v_i,v_j\}\notin\H$ for all $i<j$, or {\em dominating}, meaning that $\{v_i,v_j\}\in\H$ for all $i<j$.
A proof can be found in \cite{CH1,Gol,MP}.

\bp{order}
A graph (i.e. $2$-hypergraph) $(V,\H)$ is separable if and only if it orderable.
\ep
We point out that the above definition of orderable graphs can be naturally
extended to $k$-hypergraphs for all $k$, but while every orderable $k$-hypergraph
is separable, for all $k\geq 3$ there are $k$-hypergraphs which are separable
but not orderable, see \cite[Proposition 3.1]{Onn}.

\vskip.2cm
The following theorem resolves Question \ref{question2} for all $k$-hypergraphs with $k=1,2$.
\bt{graph}
For $k=1,2$, a $k$-hypergraph $(V,\H)$ is equatable if and only if it is exchangeable\break
(there are $E_1,E_2\in\H$ and $v_1\in E_1\setminus E_2$, $v_2\in E_2\setminus E_1$ with
$E_1\setminus\{v_1\}\cup\{v_2\}$, $E_2\setminus\{v_2\}\cup\{v_1\}\in\HC$).
\et
\bpr
If there are such $E_1,E_2,v_1,v_2$ then the hypergraph is equatable by Lemma \ref{sufficient}.
So we need to prove that if $(V,\H)$ is equatable then there are such suitable $E_1,E_2,v_1,v_2$.
For $k=1$ every $1$-hypergraph is separable with the labeling $x$ defined by $x(v):=0$
whenever $\{v\}\in\H$ and $x(v):=-1$ whenever $\{v\}\notin\H$, so this statement trivially holds.

Consider $k=2$. Since we assume that the graph is equatable, by Lemma \ref{separable_equatable}
it is not separable, which is equivalent to the graph not being threshold. Then, by a well known
characterization of threshold graph, see \cite{CH1,Gol,MP}, there are four vertices $a,b,c,d\in V$ such
that $ab,cd\in\H$ but $ac,bd\notin\H$. Taking $v_1:=a$, $v_2:=d$, $E_1:=ab$, $E_2:=cd$, we are done.
\epr

\section{Multipartite hypergraphs}

We write $\uplus_{i=1}^k V_i$ for the disjoint union of sets $V_1,\dots,V_k$,
indicating in particular that they are pairwise disjoint.
A {\em $k$-partite hypergraph} is a $k$-hypergraph $(\uplus_{i=1}^k V_i,\H)$
with a specified $k$-partition of its vertex set $V:=\uplus_{i=1}^k V_i$ such that
$\H\subseteq \{E\subseteq V: |E\cap V_i|=1,\ i=1,\dots,k\}$.

We now resolve Question \ref{question2} for $k$-partite hypergraphs for all positive $k$.
\bt{multipartite}
A $k$-partite hypergraph $(\uplus_{i=1}^k V_i,\H)$ is equatable if and only if it is exchangeable
(there are $E_1,E_2\in\H$ and $v_1\in E_1\setminus E_2$, $v_2\in E_2\setminus E_1$ with
$E_1\setminus\{v_1\}\cup\{v_2\}$, $E_2\setminus\{v_2\}\cup\{v_1\}\in\HC$).
\et
\bpr
If there are such $E_1,E_2,v_1,v_2$ then the hypergraph is equatable by Lemma \ref{sufficient}.
So we need to prove that if $(V,\H)$ is equatable then there are suitable $E_1,E_2,v_1,v_2$.

If $|\H|\leq 1$ then $(V,\H)$ is separable hence not equatable by
Lemma \ref{separable_equatable}, trivially proving the claim.
So we may assume $|\H|\geq 2$. Also, by Theorem \ref{graph}, we may assume $k\geq 3$.

Suppose for a contradiction that $|E_1\cap E_2|=k-1$ for all distinct $E_1,E_2\in\H$.
Pick any two edges $E_1\neq E_2$. Relabeling the $V_i$ if necessary, we may assume that for some
$u_1,u_2\in V_1$ and some $v_i\in V_i$ for $2\leq i\leq k$ we have $E_j=\{u_j,v_2,\dots,v_k\}$
for $j=1,2$. Then we claim that for some $U\subseteq V_1$ we have in fact
$\H=\{\{u,v_2,\dots,v_k\}:u\in U\}$. Suppose $E$ is any edge other than $E_1,E_2$.
If for some $2\leq i\leq k$ we have $v_i\notin E$ then for $j=1,2$ we have that
$|E\cap E_j|=k-1$ implies $u_j\in E$, so $u_1,u_2\in E$ contradicting $|E\cap V_1|=1$.
So the claim follows. But then $(V,\H)$ is seprabale, contradicting it being equatable,
as the following labeling shows:
$$x(v)\ :=\ \left\{
  \begin{array}{ll}
   1, & v=v_2,\dots,v_k; \\
    -(k-1), & v\in U; \\
    -k, & \hbox{otherwise}.
  \end{array}
\right.$$
So there are edges $E_1,E_2\in\H$ with $|E_1\cap E_2|\leq k-2$
and hence there are $1\leq i<j\leq k$ and
$$v_{1,i}\in (E_1\setminus E_2)\cap V_i,\ \ v_{1,j}\in (E_1\setminus E_2)\cap V_j,\ \
v_{2,i}\in (E_2\setminus E_1)\cap V_i,\ \ v_{2,j}\in (E_2\setminus E_1)\cap V_j\ .$$
Define $F_1:=E_1\setminus\{v_{1,i}\}\cup\{v_{2,j}\}$ and
$F_2:=E_2\setminus\{v_{2,j}\}\cup\{v_{1,i}\}$. Then $|F_1\cap V_j|=2$ and $|F_2\cap V_i|=2$ so
$F_1,F_2\in\HC$. Defining $v_1:=v_{1,i}$, $v_2:=v_{2,j}$ we get the desired $E_1,E_2,v_1,v_2$.
\epr

\section{Matroids}

We now consider hypergraphs which are (sets of bases of) matroids. See \cite{Oxl} for a
detailed development of the theory of matroids. A {\em $k$-matroid} is a $k$-hypergraph
$(V,\H)$ such that $\H\neq\emptyset$ and for every $E_1,E_2\in\H$ and $v_1\in E_1\setminus E_2$
there is a $v_2\in E_2\setminus E_1$ such that $E_1\setminus\{v_1\}\cup\{v_2\}\in\H$
(compare this with the second condition in Question \ref{question2}). Thus, in standard terminology,
$\H$ is the set of {\em bases} of a matroid of rank $k$ on $V$. For instance, both separable
and equatable $3$-hypergraphs demonstrated in the introduction are $3$-matroids.\break

\vskip-0.4cm
A subset $I\subseteq V$ is {\em independent} in the matroid if $I\subseteq E$ for some basis $E\in\H$.
If $I$ is any independent set and $E\in\H$ is any basis then $I$ can be {\em augmented} from $E$
to a basis, that is, $I\uplus J\in\H$ for some $J\subseteq E\setminus I$. A {\em loop} in a
$k$-matroid $(V,\H)$ is a vertex $v\in V$ that is not contained in any $E\in\H$.
A {\em coloop} is a vertex $v\in V$ that is contained in every $E\in\H$. If $v$ is not a coloop then let
$\H\backslash v:=\{E:v\notin E\in\H\}$, and if $v$ is not a loop then let
$\H/v:=\{E\setminus\{v\}:v\in E\in\H\}$. If $v$ is a coloop (hence not a loop) then let
$\H\backslash v:=\H/v$, and if $v$ is a loop (hence not a coloop) then let $\H/v:=\H\backslash v=\H$.
The {\em deletion of $(V,\H)$ by $v$} is $(V,\H)\backslash v:=(V\setminus\{v\},\H\backslash v)$,
and is a $k$-matroid if $v$ is not a coloop and a\break $(k-1)$-matroid if $v$ is a coloop.
The {\em contraction of $(V,\H)$ by $v$} is $(V,\H)/v:=(V\setminus\{v\},\H/v)$,
and is a $(k-1)$-matroid if $v$ is not a loop and a $k$-matroid if $v$ is a loop.

\subsection{Paving matroids}

A {\em paving} $k$-matroid is a $k$-matroid with the property that every $(k-1)$-subset of $V$
is independent, that is, contained in some $E\in\H$. This is a broad and fundamental class of
matroids, and in fact it is conjectured that almost all matroids are paving \cite{MNWW}. It is easy
to verify that any deletion and any contraction of a paving matroid is again a paving matroid. For
instance, the $3$-matroid $(V,\H)$ with $V=[4]$ and $\H={V\choose 3}$ is paving and separable with
the zero labeling $x$, and the $3$-matroid $(V,\H)$ with the following data is paving and equatable,
$$V=[5],\quad\H=\{123,124,134,135,145,234,235,245\},\quad
y(G)=
\left\{
  \begin{array}{rl}
    1, & G=125,135,245,345; \\
    0, & \hbox{otherwise}.
  \end{array}
\right.
$$

We now resolve Question \ref{question2} for paving $k$-matroids for all positive $k$.
\bt{paving_matroids}
A paving $k$-matroid $(V,\H)$ is equatable if and only if it is exchangeable\break
(there are $E_1,E_2\in\H$ and $v_1\in E_1\setminus E_2$, $v_2\in E_2\setminus E_1$ with
$E_1\setminus\{v_1\}\cup\{v_2\}$, $E_2\setminus\{v_2\}\cup\{v_1\}\in\HC$).
\et
\bpr
If there are such $E_1,E_2,v_1,v_2$ then the matroid is equatable by Lemma \ref{sufficient}.
So we need to prove that if $(V,\H)$ is equatable then there are suitable $E_1,E_2,v_1,v_2$.

We use induction on $k$. For $k=1,2$ this follows from Theorem \ref{graph} providing the base.
Consider $k\geq 3$. We prove the claim for $k$ by induction on $n=|V|>k$. If $n=k+1$ or $n=k+2$
then the dual $(V,\H^*)$ is a $1$-hypergraph or a $2$-hypergraph and the statement
follows from Theorem \ref{graph} and Lemma \ref{complement_dual}. Consider $n\geq k+3$.

Suppose first that there is some $v\in V$ which is not contained in any $F\in\HC$.
Then $v$ is not a coloop else $\H=\{E\in{V\choose k}:v\in E\}$ so the hypergraph
is separable with the labeling $x(v):=k-1$ and $x(u):=-1$ for all $u\neq v$, contradicting
it being equatable. So $\H\backslash v=\{E:v\notin E\in\H\}$. Since $(V,\H)$ is equatable
there is a labeling $y$ of $V\choose k$ satisfying \eqref{equatable_definition}.
This condition at $v$ implies $y(G)=0$ for all $G\in{V\choose k}$ with $v\in G$.
Therefore the restriction of $y$ to ${V\setminus\{v\}\choose k}$ shows that $(V,\H)\backslash v$
is an equatable paving $k$-matroid. Then by induction on $n$ there are suitable
$E_1,E_2\subseteq V\setminus\{v\}$ and $v_1,v_2\in V\setminus\{v\}$.
These are suitable for $(V,\H)$ as well, and the induction on $n$ follows.
So we may assume every $v\in V$ is contained in some $F\in\HC$.

Now pick any $v\in V$. As $(V,\H)$ is paving, $v$ is not a loop, so
$\H/ v=\{E\setminus\{v\}:v\in E\in\H\}$. Now, if $(V,\H)/v$ is equatable then,
since it is a paving $(k-1)$-matroid, by induction on $k$ there are suitable
$E_1,E_2\subseteq V\setminus\{v\}$ and $v_1,v_2\in V\setminus\{v\}$. Then
$E_1\cup\{v\},E_2\cup\{v\},v_1,v_2$ are suitable for $(V,\H)$,
and the induction on $n$ follows. Assume then that $(V,\H)/v$ is not equatable and hence by
Lemma \ref{separable_equatable} is separable. Let $x$ be a suitable labeling of $V\setminus\{v\}$
and index the vertices in $V\setminus\{v\}$ so that $x(v_1)\leq\cdots\leq x(v_{n-1})$.
If $x(v_1)+\cdots+x(v_{k-2})+x(v_{n-1})<0$ then $x(v_1)+\cdots+x(v_{k-2})+x(v_i)<0$
for all $k-1\leq i\leq n-1$. Then $v_1\cdots v_{k-2}v_i\notin\H/v$ and hence
$v_1\cdots v_{k-2}v_iv\notin\H$ for all $i$. But $v_1\cdots v_{k-2}v$ is independent in $(V,\H)$
since it is paving, so must be contained in some basis, which is a contradiction.
So $x(v_1)+\cdots+x(v_{k-2})+x(v_{n-1})\geq 0$ and hence for all
$1\leq i_1<\cdots<i_{k-2}\leq n-2$ we have that
$x(v_{i_1})+\cdots+x(v_{i_{k-2}})+x(v_{n-1})\geq 0$ and hence
$v_{i_1}\cdots v_{i_{k-2}}v_{n-1}\in\H/v$ and $v_{i_1}\cdots v_{i_{k-2}}v_{n-1}v\in\H$. Letting
$u:=v_{n-1}$, this means that any $k$-subset of $V$ which contains $u,v$ is in $\H$.

Now, as proved above, there are $F_1,F_2\in\HC$ with $u\in F_1$ and $v\in F_2$.
Since $(V,\H)$ is paving, $F_1\setminus\{u\}$ is independent in $\H$. Since $F_1\cup\{v\}$
contains $u,v$, by what just proved it contains a basis in $\H$ from which $F_1\setminus\{u\}$
can be augmented\ to a basis. This implies that $E_1:=F_1\setminus\{u\}\cup\{v\}\in\H$.
A similar argument shows that $E_2:=F_2\setminus\{v\}\cup\{u\}\in\H$. Defining $v_1:=v$,
$v_2:=u$ we get the desired $E_1,E_2,v_1,v_2$ and the inductions on $n$ and $k$ follow.
\epr

\subsection{Three-matroids}

We need two lemmas and some more matroid terminology as follows.
\bl{loops}
A $k$-matroid $(V,\H)$ with a loop $v\in V$ is separable if and only if $(V,\H)\backslash v$ is.
\el
\bpr
Since $v$ is a loop we have $\H\backslash v=\H$. In one direction, suppose $(V,\H)$
is separable with labeling $x:V\rightarrow\R$. Then its restriction to $V\setminus\{v\}$ shows
$(V,\H)\backslash v$ is separable, since a $k$-subset $E\subseteq V\setminus\{v\}$
satisfies $x(E)\geq 0$ if and only if $E\in\H=\H\backslash v$. In the other direction,
suppose $(V,\H)\backslash v$ is separable with labeling $x:V\setminus\{v\}\rightarrow\R$.
Extend it to $V$ by setting $x(v):=-r$ where $r$ is a sufficiently large positive number.
Then $(V,\H)$ is separable since
$$\{E\subseteq V:|E|=k,\ x(E)\geq 0\}\ =\
\{E\subseteq V\setminus\{v\}:|E|=k,\ x(E)\geq 0\}\ =\ \H\backslash v\ =\ \H\ .\mepr $$

\vskip.2cm
Let $(V,\H)$ be a matroid with no loops. A {\em line} of the matroid is a subset
$\emptyset\neq L\subseteq V$ such that $uv$ is not independent for all distinct $u,v\in L$ and $uv$
is independent for all $u\in L$ and $v\in V\setminus L$. The line is
{\em nontrivial} if $|L|\geq 2$. The set $V$ equals the disjoint union of the lines.

\bl{lines}
A matroid $(V,\H)$ with no loops and at least two nontrivial lines is exchangeable
(there are $E_1,E_2\in\H$ and $v_1\in E_1\setminus E_2$, $v_2\in E_2\setminus E_1$
with $E_1\setminus\{v_1\}\cup\{v_2\}$, $E_2\setminus\{v_2\}\cup\{v_1\}\in\HC$).
\el
\bpr
Let $(V,\H)$ be a $k$-matroid with no loops and two distinct nontrivial lines $L_1,L_2$.
Pick any distinct $u_1,v_2\in L_1$ and $u_2,v_1\in L_2$. Since $u_1v_1$ is independent
it is contained in some basis $E_1\in\H$. Since $u_2v_2$ is independent it is contained in
some basis $E_2\in\H$. Now $E_1\setminus\{v_1\}\cup\{v_2\}\in\HC$ since it contains $u_1v_2$
which is not independent. Similarly, $E_2\setminus\{v_2\}\cup\{v_1\}\in\HC$ since it
contains $u_2v_1$ which is not independent. The lemma follows.
\epr

We now resolve Question \ref{question2} for all $3$-matroids.
\bt{three-matroids}
A $3$-matroid $(W,\H)$ is equatable if and only if it is exchangeable\break
(there are $E_1,E_2\in\H$ and $v_1\in E_1\setminus E_2$, $v_2\in E_2\setminus E_1$ with
$E_1\setminus\{v_1\}\cup\{v_2\}$, $E_2\setminus\{v_2\}\cup\{v_1\}\in\HC$).
\et
\bpr
If there are such $E_1,E_2,v_1,v_2$ then the matroid is equatable by Lemma \ref{sufficient}.
So we need to prove that if $(W,\H)$ is equatable then there are suitable $E_1,E_2,v_1,v_2$.

Deleting all loops if any one after the other we obtain a loopless $3$-matroid $(V,\H)$
with $n:=|V|\geq 3$, which is equatable if and only if $(W,\H)$ is by Lemma \ref{loops}.
And, if $E_1,E_2$ and $v_1,v_2$ are good for $(V,\H)$, then they are also good for $(W,\H)$.
So it suffices to prove the claim for $(V,\H)$. If $n=3$ then $\H=\{V\}$ so $(V,\H)$
is separable. If $n=4$ or $n=5$ then the dual $(V,\H^*)$ is a $1$-hypergraph or a
$2$-hypergraph so the claim follows from Theorem \ref{graph}.

Consider then $n\geq 6$. If $(V,\H)$ has no nontrivial lines then every
$2$-subset of $V$ is independent, and therefore the matroid is paving and the
claim follows from Theorem \ref{paving_matroids}. If the matroid has at least two
distinct nontrivial lines then the claim follows by Lemma \ref{lines}.

So assume there is exactly one nontrivial line $L$. Let $U:=V\setminus L$.
Pick any $v\in L$. Since $v$ is not a loop, we have $\G:=\H/v=\{ab:abv\in\H\}$.
Consider the graph $(U,\G)$. If it is equatable then by Theorem \ref{graph} there are
suitable $E_1,E_2\subseteq U$ and $v_1,v_2\in U$ and then $E_1\cup\{v\},E_2\cup\{v\}$ and
$v_1,v_2$ are suitable for $\H$ and we are done. So assume $(U,\G)$ is separable.

Consider any three distinct vertices $a,b,c\in U$ (if any) and let $i:=|\{ab,ac,bc\}\cap\G|$ be the
number of edges among $ab,ac,bc$ contained in $\G$. Suppose $i=0$ and suppose for a contradiction
that $abc\in\H$. Then we must be able to augment the non loop $v$ from $abc$ to a basis.
But $ab,ac,bc\notin\G$ so $abv,acv,bcv\notin\H$, a contradiction. So if $i=0$ then $abc\notin\H$.
Next suppose for a contradiction that $i=1$, say $ab\in\G$, which implies $abv\in\H$.
Since $v\in L$ and $c\notin L$ we have that $cv$ is independent and so we must be able to
augment it from $abv$ to a basis. So we must have either $acv\in\H$ or $bcv\in\H$ which
implies either $ac\in\G$ or $bc\in\G$ contradicting $i=1$. So $i\neq 1$. Now suppose $i\geq 2$,
say $ab,bc\in\G$ which implies $abv,bcv\in\H$. Pick any $w\in L\setminus\{v\}$.
Since $bw$ is independent we must be able to augment it to a basis from $bcv$, but
$vw$ is not independent so we find that this basis must be $bcw\in\H$. Now, if $abc\notin\H$
then taking $E_1:=abv$, $E_2:=bcw$, $v_1:=v$, $v_2:=c$ the claim follows again.

We claim that indeed there must be some $a,b,c\in U$ with $i\geq 2$ and $abc\notin\H$ so that
we are done. Assume, for a contradiction, that this is not the case. Then, under this assumption,
the following hold: if $i=0$ then $abc\notin\H$; always $i\neq 1$; and $abc\in\H$ whenever $i=2,3$.

Now, for any $a\in U$ and any distinct $w,z\in L$ we have that $wz$ is not independent so
$awz\notin\H$. Also, for every distinct $a,b\in U$ and $w\in L$, we have that $abw\in\H$ if and only
if $abv\in\H$ if and only if $ab\in\G$. In particular, since $\H$ is a matroid and hence nonempty,
we see that $\G$ cannot be empty, so has at least one edge and $(U,\G)$ has $r:=|U|\geq 2$ vertices.

Now, since $(U,\G)$ is separable, by Lemma \ref{order} it has an ordering
$u_1,\dots,u_r$ of $U$ such that each $u_i$ is either isolated or dominating, with $u_1$ declared
isolated. Now, if there are $1<i<j\leq r$ with $u_i$ dominating and $u_j$ isolated then
$\{u_1u_i,u_1u_j,u_iu_j\}\cap\G=\{u_1u_i\}$ which is impossible as shown above. So
for some $1\leq i<r$ we have that $u_j$ is isolated for $j\leq i$ and dominating for $j>i$.
Let $U_0:=\{u_j:j\leq i\}$ and $U_1:=\{u_j:j>i\}$ so that $U=U_0\uplus U_1$.

We claim that $(V,\H)$ is separable, contradicting it being equatable. Consider the labeling
$$x(u)\ :=\
\left\{
  \begin{array}{rl}
    -1, & u\in U_0; \\
    3, & u\in U_1; \\
    -2, & u\in L.
  \end{array}
\right.
$$
Consider any $3$-subset $S\subset V$. We show that $S\in\H$ if and only if $x(S)\geq 0$.
Suppose first that $|S\cap L|=0$. If $|S\cap U_1|=0$ then $|{S\choose 2}\cap\G|=0$ so $S\in\HC$
and indeed $x(S)=-3<0$. If $|S\cap U_1|\geq 1$ then $|{S\choose 2}\cap\G|\geq 2$ so
$S\in\H$ and indeed $x(S)\geq 3-1-1\geq 0$. Next suppose that $|S\cap L|=1$.
If $S\cap U\in\G$ then $S\in\H$ and indeed $|S\cap U_1|\geq 1$ hence $x(S)\geq 3-1-2\geq 0$.
If $S\cap U\notin\G$ then $S\in\HC$ and indeed $|S\cap U_1|=0$ hence $x(S)=-2-1-1<0$.
Finally suppose that $|S\cap L|\geq 2$. Then $S\in\HC$ and indeed $x(S)\leq -2-2+3<0$.

We conclude there are $a,b,c\in U$ with $i:=|\{ab,ac,bc\}\cap\G|\geq 2$ and $abc\notin\H$, so there
are suitable $E_1,E_2$ and $v_1,v_2$ as explained above. This completes the proof.
\epr

\subsection{Binary matroids}

We need some more matroid terminology. A {\em circuit} in a matroid $(V,\H)$ is a subset
$C\subseteq V$ that is not independent but all its proper subsets are independent.
If $C_1,C_2$ are two circuits with $v\in C_1\cap C_2$ and $u\in C_1\setminus C_2$ then
there is another circuit $C$ such that $u\in C\subseteq (C_1\cup C_2)\setminus\{v\}$.
If $E\in\H$ and $v\in V\setminus E$ then there is a unique circuit $C(E,v)$ satisfying
$C(E,v)\subseteq E\uplus\{v\}$. A {\em binary $k$-matroid} is a $k$-matroid with the
property that for any two distinct circuits $C_1,C_2$, their {\em symmetric difference}
$C_1\Delta C_2:=(C_1\setminus C_2)\uplus(C_2\setminus C_1)$ is a disjoint union of circuits. This
is a broad and well studied class that includes in particular all {\em graphic matroids}, in which
$V$ is the set of edges of a graph and $\H$ is the set of maximal forests in the graph, see \cite{Oxl}.

We now resolve Question \ref{question2} for binary $k$-matroids for all positive $k$.
\bt{binary}\hskip-.1cm
A binary $k$-matroid $(W,\H)$ is equatable if and only if it is exchangeable\break
(there are $E_1,E_2\in\H$ and $v_1\in E_1\setminus E_2$, $v_2\in E_2\setminus E_1$
with $E_1\setminus\{v_1\}\cup\{v_2\}$, $E_2\setminus\{v_2\}\cup\{v_1\}\in\HC$).
\et
\bpr
If there are such $E_1,E_2,v_1,v_2$ then the matroid is equatable by Lemma \ref{sufficient}.
So we need to prove that if $(W,\H)$ is equatable then there are suitable $E_1,E_2,v_1,v_2$.

Deleting all loops if any one after the other we obtain a loopless $k$-matroid $(V,\H)$
and $n:=|V|\geq k$, which is equatable if and only if $(W,\H)$ is by Lemma \ref{loops}.
And, if $E_1,E_2,v_1,v_2$ are good for $(V,\H)$, then they are also good for $(W,\H)$.
So it suffices to prove the claim for $(V,\H)$. If $n=k$ then $\H=\{V\}$ so $(V,\H)$ is
separable. If $n=k+1$ then the dual $(V,\H^*)$ is a $1$-hypergraph so separable,
see proof of Theorem \ref{graph}, and hence so is $(V,\H)$ by Lemma \ref{complement_dual}.

Consider $n\geq k+2$. If $(V,\H)$ has at least two distinct nontrivial lines then by
Lemma \ref{lines} there are suitable $E_1,E_2,v_1,v_2$ and we are done. Suppose then
there is at most one nontrivial line and let $L$ be one with maximum $l:=|L|\geq 1$.
Since $|E\cap L|\leq 1$ for all $E\in\H$ we have $n\geq l+k-1$. If $n=l+k-1$ then
$\H=\{E\in{V\choose k}:|E\cap L|=1\}$ so the labeling $x(u):=1$ for $u\in V\setminus L$
and $x(u):=-(k-1)$ for $u\in L$ shows $(V,\H)$ is separable. So assume $n\geq l+k$. Pick any
$E\in\H$ with $|E\cap L|=1$. Then there is some $v_2\in V\setminus(E\cup L)$ and $\{v_2\}$
must be a trivial line. Pick any $v_1\in V\setminus(E\uplus\{v_2\})$. For $i=1,2$ let $C_i:=C(E,v_i)$
and $P_i:=C_i\setminus\{v_i\}$. We claim $P_1\Delta P_2\neq\emptyset$. Indeed, otherwise we
get a contradiction, since then $C_1\Delta C_2=\{v_1,v_2\}$; but $C_1\Delta C_2$ is a
disjoint union of circuits since the matroid is binary, while $v_1v_2$ is independent.
So the claim is true, and we can assume, say, that there is some element $u_1\in P_1\setminus P_2$.

\vskip.1cm
Suppose there is also an element $u_2\in P_2\setminus P_1$. Define the following sets,
$$E_1:=E\uplus\{v_1\}\setminus\{u_1\},\ \ E_2:=E\uplus\{v_2\}\setminus\{u_2\},\ \
F_1:=E_1\setminus\{v_1\}\uplus\{v_2\},\ \ F_2:=E_2\setminus\{v_2\}\uplus\{v_1\}\ .$$
Then, for $i=1,2$, we have that $E_i\in\H$ since $u_i\in C_i$, and $F_i\in\HC$ since
$C_{3-i}\subseteq F_i$.

Now suppose $P_2\subsetneq P_1$. Since $v_2v$ is independent for all $v\neq v_2$ we have
$|C_2|\geq 3$ and hence $|P_2|\geq 2$. Pick any two distinct vertices $u_2,w_2\in P_2$.
Let $G:=E\uplus\{v_1\}\setminus\{u_1\}\in\H$. Since $C_2\subseteq G\uplus\{v_2\}$ we find that
$C(G,v_2)=C_2$ and therefore
$$E'\ :=\ G\uplus\{v_2\}\setminus\{u_2\}\ =\ E\uplus\{v_1,v_2\}\setminus\{u_1,u_2\}\in\H\ .$$
For $i=1,2$ let $C'_i:=C(E',u_i)$ and $P_i':=C'_i\setminus\{u_i\}$.
As we have seen, $C'_2=C_2$. Now consider $C_1\Delta C_2=(P_1\setminus P_2)\uplus\{v_1,v_2\}$
which must be a disjoint union of circuits. Since $u_2\in C_1\cap C_2$ and
$u_2\notin(P_1\setminus P_2)\uplus\{v_i\}\subset E\uplus\{v_i\}$ we conclude that
$(P_1\setminus P_2)\uplus\{v_i\}$ contains no circuit for $i=1,2$. So we must have
$(P_1\setminus P_2)\uplus\{v_1,v_2\}=C'_1$. So $v_1\in P'_1\setminus P'_2$ and
$w_2\in P'_2\setminus P'_1$. Define
$$E'_1:=E'\uplus\{u_1\}\setminus\{v_1\},\ \ E'_2:=E'\uplus\{u_2\}\setminus\{w_2\},\ \
F'_1:=E'_1\setminus\{u_1\}\uplus\{u_2\},\ \ F'_2:=E'_2\setminus\{u_2\}\uplus\{u_1\}\ .$$
Then $E'_1,E'_2\in\H$ since $v_1\in C'_1$ and $w_2\in C'_2$, and for $i=1,2$ we have
$F'_i\in\HC$ since $C'_{3-i}\subseteq F'_i$.

\vskip.1cm
So in either case there are two suitable bases and non bases and the claim follows.
\epr

\section{An equatable non exchangeable hypergraph}

As mentioned in the introduction, the conditions in Questions \ref{question1} and \ref{question2}
do not hold for all hypergraphs. The following example from \cite{RRST} gives a non separable
$3$-hypergraph which does not satisfy these conditions. Using the notion of equatability and
the assertion of Lemma \ref{separable_equatable} that a hypergraph is equatable if and only
if it is non separable, we provide a simpler demonstration that this example is indeed valid.
Let $(V,\H)$ be the following $3$-hypergraph,
$$V=[9],\quad\H=\{ij9\ :\ i<j<9,\ 1\leq i,\ 4\leq j\}
\cup\{ijk\ :\ i<j<k,\ 2\leq i,\ 5\leq j,\ 7\leq k\}\ .$$
It is easy to verify that the following labeling $y:{[9]\choose3}\rightarrow\{0,1\}$
shows that it is equatable,
$$
y(G)=
\left\{
  \begin{array}{rl}
    1, & G=149,178,239,267,358,456; \\
    0, & \hbox{otherwise}.
  \end{array}
\right.
$$
Now, consider any $E_1,E_2\in\H$ and any $v_1\in E_1\setminus E_2$,
$v_2\in E_2\setminus E_1$. Reindexing if necessary we may assume $v_2>v_1$.
But then, by the definition of $\H$, we must have that $E_1\setminus\{v_1\}\cup\{v_2\}\in\H$.
So $(V,\H)$ is not exchangeable and the conditions of Questions \ref{question1}
and \ref{question2} fail. Note that $(V,\H)$ is not a matroid: both $E_1:=149$, $E_2:=257$
are in $\H$ but there is no $v_2\in E_2$ such that $E_1\setminus\{9\}\cup\{v_2\}\in\H$.
So the condition of Question \ref{question2} might yet hold for all matroids. Also, $y$ is
$0,1$ valued with $6=2k$ non zeros, so Question \ref{question3} might yet have a positive answer.

\section{Monotone hypergraphs}

As mentioned in the introduction, related to Question \ref{question2} is the notion of monotonicity.
In the context of Boolean functions this had been studied extensively, see \cite[Section 8.8]{CH2}.
Here we consider the uniform analog of this property, relevant in the context of $k$-hypergraphs.

Let $(V,\H)$ be a $k$-hypergraph. For $R_1,R_2\subseteq V$ with $|R_1|=|R_2|$ put $R_1\leq R_2$
if $S\cup R_1\in\H$ implies $S\cup R_2\in\H$ for all $S\subseteq V\setminus(R_1\cup R_2)$.
The hypergraph is (uniformly) {\em $r$-monotone} if for all $R_1,R_2\subseteq V$ such that
$|R_1|=|R_2|$ and $|R_1\cup R_2|\leq r$, $\mbox{either}\ R_1\leq R_2\ \mbox{or}\ R_2\leq R_1$.

Every separable $k$-hypergraph is $r$-monotone for all $k$ and $r$. To see this,
suppose $(V,\H)$  is separable with a suitable labeling $x:V\rightarrow\R$.
Consider any $R_1,R_2\subseteq V$ with $|R_1|=|R_2|$. Reindexing if necessary, we may assume
$x(R_1)\leq x(R_2)$. Then for all $S\subseteq V\setminus(R_1\cup R_2)$ we have that if
$S\cup R_1\in\H$ then  $x(S\cup R_2)\geq x(S\cup R_1)\geq 0$ so $S\cup R_2\in\H$. So $R_1\leq R_2$.

\bp{monotone} A $k$-hypergraph is $2$-monotone if and only if it is {\em not} exchangeable
(for all $E_1,E_2\in\H$ and $v_1\in E_1\setminus E_2$, $v_2\in E_2\setminus E_1$,
either $E_1\setminus\{v_1\}\cup\{v_2\}\in\H$ or $E_2\setminus\{v_2\}\cup\{v_1\}\in\H$).
\ep
\bpr
Suppose first the hypergraph is $2$-monotone. Consider any relevant $E_1,E_2$ and $v_1,v_2$.
Reindexing if necessary, we may assume $\{v_1\}\leq\{v_2\}$. Let $S:=E_1\setminus\{v_1\}$.
Then $S\cup\{v_1\}=E_1\in\H$ implies $E_1\setminus\{v_1\}\cup\{v_2\}=S\cup\{v_2\}\in\H$.
So the hypergraph is not exchangeable.

Conversely, suppose the hypergraph is not $2$-monotone. Then there must be distinct
singletons $R_1=\{v_1\}$, $R_2=\{v_2\}$, and $S_1,S_2\subseteq V\setminus(R_1\cup R_2)$,
such that $S_1\cup R_1,S_2\cup R_2\in\H$ but $S_1\cup R_2,S_2\cup R_1\in\HC$.
Then the hypergraph is exchangeable since we have
$$E_1\ :=\ S_1\cup R_1\ \in\ \H,\quad E_2\ :=\ S_2\cup R_2\ \in\ \H\ ,$$
$$E_1\setminus\{v_1\}\cup\{v_2\}\ =\ S_1\cup R_2\ \in\ \HC,\quad
E_2\setminus\{v_2\}\cup\{v_1\}\ =\ S_2\cup R_1\ \in\ \HC\ .
\mepr
$$

\vskip.2cm
Since every separable hypergraph is $r$-monotone for all $r$, in particular every
separable matroid is $2$-monotone. But which $2$-monotone matroids are separable?
Our results in Theorem \ref{main} combined with Lemma \ref{separable_equatable}
and Proposition \ref{monotone} show that $2$-monotone paving matroids, $2$-monotone binary
matroids, and $2$-monotone $3$-matroids, are all separable. And in \cite{GK} it is shown
that $2$-monotone matroids which are moreover {\em non uniform} $3$-monotone are separable.
Could these results be unified and strengthened to positively answer the following question:

\bq{monotone_question}
Is it true that any matroid is separable if and only if it is $2$-monotone?\break
Equivalently, is it true that any matroid is equatable if and only if it is exchangeable?\break
\eq

\section{Remarks on complexity}

We conclude with some remarks on the complexity of deciding if a hypergraph is separable.

First (see also \cite{RRST}), for any fixed $k$, the number $O(n^k)$ of $k$-subsets of $V$ is
polynomial in $n=|V|$. So, given the list of edges in $\H$, it can be decided in polynomial time if
$(V,\H)$ is separable by writing the system of linear inequalities \eqref{primal} in
Lemma \ref{separable_equatable} and deciding if it has a solution using
linear programming, which can be done in polynomial time, see e.g. \cite{Sch}.

Second, assume $k$ is a {\em variable} part of the input. Let $m:=|\H|$ and assume the hypergraph
is given by a list of the $m$ edges in $\H$. The systems \eqref{primal}, \eqref{dual} are then
of exponential size $\Omega\left({n\choose k}\right)$ and cannot be used to directly decide if
the hypergraph is separable. But we can check in time polynomial in $m$, which may be much
smaller than $n\choose k$, if it is exchangeable, by checking for each of the $O(m^2)$
pairs of edges $E_1,E_2\in\H$ and each of the $O(k^2)$ vertices $v_1\in E_1\setminus E_2$
and $v_2\in E_2\setminus E_1$ whether or not $E_1\setminus\{v_1\}\cup\{v_2\}$,
$E_2\setminus\{v_2\}\cup\{v_1\}\in\HC$. So if $\H$ is in a class satisfying the condition of
Question \ref{question2}, we can efficiently check if it is exchangeable and hence equatable,
and hence efficiently decide indirectly, by Lemma \ref{separable_equatable}, if it is separable.
We note that \cite{RRST} describes an algorithm for deciding if a hypergraph is $2$-monotone
(therein $2$-monotonicity is called {\em property $T_3$}) with exponential running time
$O(n^k)$. Since by Proposition \ref{monotone} a hypergraph is $2$-monotone if and only
if it is not exchangeable, the above procedure gives a simpler, polynomial rather than exponential,
algorithm for deciding $2$-monotonicity.

Third, assume $k$ is again a variable part of the input, and the hypergraph, which may have
a number of edges which is exponential in $k$ and $n$, is a matroid which is presented implicitly
by an {\em independence oracle} that, queried on any subset $I\subseteq V$, replies YES if $I$
is independent and NO if it is not independent. We claim that, even if the matroid is paving,
then exponentially many queries are needed in general to decide if it is separable, but if the
matroid is binary, then it can be decided if it is separable using polynomially many queries.

Consider paving $k$-matroids with $n=2k$ and $k\geq 2$. Suppose an algorithm
trying to solve the decision problem makes less than $2^k-1\leq{1\over 2}{2k\choose k}$ queries.
Then there is a pair of disjoint $k$-sets $F_1,F_2$ about which the algorithm did not query
the oracle. Let $E_1,E_2\neq F_1,F_2$ be two other disjoint $k$-sets. Consider the complete
$k$-hypergraph $(V,\H_1)$ with $\H_1:={V\choose k}$ and the $k$-hypergraph $(V,\H_2)$ with
$\H_2:=\H_1\setminus\{F_1,F_2\}$. Then $(V,\H_1)$ is separable with the identically zero
labeling $x$, whereas $(V,\H_2)$ is equatable by Lemma \ref{sufficient} since $E_1,E_2\in\H_2$,
$F_1,F_2\in\HC_2$, $E_1\cap E_2=\emptyset=F_1\cap F_2$, and $E_1\cup E_2=V=F_1\cup F_2$.
Clearly $(V,\H_1)$ is a paving matroid and we claim that so is $(V,\H_2)$.
It suffices to show that for any $(k-1)$-set $I\subset V$ and any $E\in\H_2$ there is some
$v\in E\setminus I$ so that $I\uplus\{v\}\in\H_2$. If $I\subset E$ this is clear. Otherwise
there are distinct $v_1,v_2\in E\setminus I$. Since $G_1:=I\uplus\{v_1\}$ and $G_2:=I\uplus\{v_2\}$
are distinct and not disjoint $k$-sets we have that $\{G_1,G_2\}\neq\{F_1,F_2\}$ so at least one
of $G_1,G_2$ is in $\H_2$. So $(V,\H_2)$ is a paving matroid. But now, whether the
oracle presents $(V,\H_1)$ or $(V,\H_2)$, it answers YES to all queries made
by the algorithm, so the algorithm cannot distinguish between these two matroids
and cannot tell whether the one presented by the oracle is separable or equatable.

Now suppose that $(W,\H)$ is a binary $k$-matroid. The proof of Theorem \ref{binary}
leads to an efficient algorithm as follows. Using $|W|$ queries we can
identify all loops if any. We let $(V,\H)$ be the matroid obtained by deleting them one
after the other and let $n:=|V|\geq k$. By Lemma \ref{loops} the original matroid is
separable if and only if the deleted matroid is. Moreover, the oracle for $(W,\H)$ is good
also for $(V,\H)$. So it suffices to solve the decision problem for the loopless $(V,\H)$.
If $n\leq k+1$ then $(V,\H)$ is separable. Assume $n\geq k+2$. By asking the oracle about
all $O(n^2)$ vertex pairs $\{u,v\}\subseteq V$ we can identify all the lines in $(V,\H)$.
If there are two nontrivial lines then $(V,\H)$ is equatable by Lemma \ref{lines}.
Otherwise let $l:=|L|$ be the largest cardinality of a line. Then, as shown in the proof
of Theorem \ref{binary}, if $n=l+k-1$ then $(V,\H)$ is separable, whereas if $n\geq l+k$
then $(V,\H)$ is equatable. So we can decide if $(V,\H)$ and hence also
$(W,\H)$ are separable using polynomially many queries to the oracle.

\section*{Acknowledgments}

D. Deza was supported by scholarships from the University of Toronto and from the Technion.
S. Onn was supported by a grant from the Israel Science Foundation and by the Dresner chair.\\


\begin{thebibliography}{}

\bibitem{BLSWZ}
Bj\"{o}rner, A., Las Vergnas, M., Sturmfels, B., White, N., Ziegler, G.M.:
Oriented matroids. Cambridge University Press (1999)

\bibitem{CH1}
Chvátal, V., Hammer, P.L.:
Aggregation of inequalities in integer programming.
Annals of Discrete Mathematics 1:145--162 (1977)

\bibitem{CH2}
Crama, Y., Hammer, P.L.:
Boolean Functions. Cambridge University Press (2011)

\bibitem{DLMO}
Deza, A., Levin, A., Meesum, S.M., Onn, S.:
Optimization over degree sequences.
SIAM Journal on Discrete Mathematics 32:2067--2079 (2018)

\bibitem{FKPT}
Frosini, A., Kocay, W.L., Palma, G., Tarsissi, L.:
On null $3$-hypergraphs.
Discrete Applied Mathematics 303:76--85 (2021).

\bibitem{GJ}
Garey, M.R., Johnson, D.S.:
Computers and Intractability. Freeman (1979)

\bibitem{GK}
Giles, R., Kannan, R.:
A characterization of threshold matroids.
Discrete Mathematics 30:181--184 (1980)

\bibitem{Gol}
Golumbic, M.C.:
Algorithmic Graph Theory and Perfect Graphs. Academic Press (1980)

\bibitem{Hax}
Haxell, P.E.:
A condition for matchability in hypergraphs.
Graphs and Combinatorics 11:245--248 (1995)

\bibitem{MP}
Mahadev, N.V.R., Peled, U.N.:
Threshold Graphs and Related Topics.
Annals of Discrete Mathematics 56, North-Holland (1995)

\bibitem{MNWW}
Mayhew, D., Newman, M., Welsh, D.J.A., Whittle, G.:
On the asymptotic proportion of connected matroids.
European Journal of Combinatorics 32:882--890 (2011)

\bibitem{Onn}
Onn, S.:
Matching orderable and separable hypergraphs.
Optimization Letters 16:1393--1401 (2022)

\bibitem{Oxl}
Oxley, J.G.:
Matroid Theory. Oxford University Press (2011)

\bibitem{RRST}
Reiterman, J., Rödl, V., Šiňajová, E., Tůma, M.:
Threshold hypergraphs. Discrete Mathematics 54:193--200 (1985)

\bibitem{Sch}
Schrijver, A.:
Theory of Linear and Integer Programming. Wiley (1986)

\end{thebibliography}
\end{document}